\documentclass[12pt]{amsart}
\usepackage{amsmath,amssymb,enumerate}
\usepackage{epsf,epsfig,amsfonts,graphicx,color}  
 \usepackage[applemac]{inputenc} 
\usepackage{hhline}
\usepackage{a4wide,epsfig}
 \numberwithin{equation}{section}   
\newtheorem{theorem}{Theorem}

\theoremstyle{definition}
\newtheorem{remark}{Remark}

\newcommand\Supp{\textrm{Supp}}
\newcommand\const{\textrm{const}}

\newcommand{\comm}[1]{}

\newcommand\op[1]{\mathop{\rm #1}\nolimits}

\newcommand\R{{\mathbb R}}

\def\R{\mathord{\mathbb R}}

\def\2{\frac{1}{2}}
\def\3{{\ss}}

\def\.{\cdot}

\def\<{\langle}
\def\>{\rangle}

\def\beq{\begin{equation}}
\def\eeq{\end{equation}}
\def\bea{\begin{eqnarray}}
\def\eea{\end{eqnarray}}
\def\bsm{\left(\begin{smallmatrix}}
\def\esm{\end{smallmatrix}\right)}
\def\bpm{\begin{pmatrix}}
\def\epm{\end{pmatrix}}

\newcommand{\weg}[1]{}

\date{}
\title[Integrable systems on  the suspensions of toric automorphim]{On integrable natural Hamiltonian systems on the suspensions of toric automorphism.}
\author{Vladimir S. Matveev}
\begin{document}
\maketitle

{\bf Abstract:} We correct  a mistake in the main statement of  \cite{liu}. \\[.1cm]

We consider $M^{n+1} = T^n \times \mathbb{R}$. We think that the torus  $T^n$ is equipped with the standard  $2\pi-$periodic 
coordinates $x_1,...,x_n$; together with the standard coordinate $x_{n+1}$ on $\mathbb{R}$ we have then a coordinate system on $M$. We take a matrix $A\in \textrm{Mat}(n,n,\mathbb{R})$ and  consider the $n\times n$-matrix 
$Q(x_{n+1})= (Q(x_{n+1})_{ij})$ given by 
$$Q(x_{n+1})= \textrm{exp}(x_{n+1} A)^t \textrm{exp}(x_{n+1} A),$$ where {\it``t''} stays for  ``transposed''.  This is a positively definite symmetric matrix. 
Next, we consider the  Riemannian metric $g$ on $M$ given by 
$$
g= \sum_{i,j=1}^nQ_{ij} dx_idx_j + dx_{n+1}^2
$$ 
and a smooth function $V = V(x_{n+1})$, which will play the role of the potential energy.
  We assume that  $V$ is periodic with  period 1, i.e. $V(x_{n+1} + 1)= V(x_{n+1})$. 
This in particular implies that the function  has a minimal   and a maximal values which we denote by $
V_{min}$ and $V_{max}$.

We consider the natural  Hamiltonian  system on $T^*M$  with   $H(p,x)= \frac{1}{2}\sum_{i,j=1}^{n+1} g^{ij}p_ip_j + V$, where $p_1,..., p_{n+1}$ are the momenta corresponding to the coordinates $x_1,...,x_{n+1}$. 
For every $h\in  \mathbb{R}$, we denote by $E_h$ the set 

$$
E_h:= \{ \xi\in T^*M\mid H(\xi)= h\}. 
$$
Each $E_h$ is invariant with respect to the Hamiltonian system. The set $E_h$ for $h<V_{min}$ is empty. For  each $h$ such that 
 $V_{min} \le h < V_{max}$,  the set $E_h$  
consists of infinitely many connected compact components. 
\begin{theorem}   
For every $h$ such that 
$V_{min} \le h<  V_{max}$ the restriction of the Hamiltonian system to any connected component of  
$E_h$  has zero topological entropy. 
\end{theorem}

  \begin{remark}  Theorem 1 evidently contradicts \cite[Theorem 1.1(c)]{liu}, from which  it follows that for almost every $h\in [V_{min}, V_{max}]$ the restriction of the Hamiltonian system to any 
    connected component of $E_h$  has positive topological entropy. The mistake in the proof of  \cite[Theorem 1.1(c)]{liu} is  hidden in the proof of  the case 2 (see pages 316-317). There, F. Liu and X. Zhang consider the invariant subset given by the condition 
     $p_1=p_2=...=p_n=0$.   They  claim that the Poincare   mapping corresponding to the restriction of the Hamiltonian system to this subset has positive  topological entropy.

     Careful analysis shows though  that this claim is wrong.  Indeed,    in this situation the Poincare map is the identical map, and its topological entropy  vanishes.    One can have an impression that Liu and Zhang thought that every  trajectory  of this subsystem\footnote{after an appropriate  reparameterization,   every   trajectory   has  locally the form  $\left(x_1=\const_1,..., x_n=\const_n, x_{n+1}=s, p_1=0, ...,p_n=0, p_{n+1} =\pm  \sqrt{2(h- V(s))} \right)$}, after reaching the point $x_{n+1}= s_2$, `jumps' to the point $x_{n+1}=s_1$, which is of cause not the case.  In fact, 
     each  trajectory returns along the same path with the reverse speed,  and  the second intersection with  the       Poincare section, if it exists, coincides with   the first one.      
  
Note that  the result of Liu and Zhang  concerns   the suspentions of toric automorphisms.  
That is,  they consider    the action of the group $\mathbb{Z}$ on $M$ generated by $$\phi(x_1,...,x_{n+1})= (\underbrace{(x_1,...,x_n)\textrm{exp}(A)^t}_{\textrm{$n$ components}} , x_{n+1}-1).$$ 
The suspension  of toric automorphism is 
 the  quotient $\widetilde M= M/_\mathbb{Z}$.   Certain additional assumptions on $A$  implying
that the metric $g$ is well defined and 
  that the topological entropy is positive for $h>V_{max}$ are assumed in \cite{liu}, 
  see \cite{liu} for details.
Since for $h\in [V_{min}, V_{max}) $ the  projection of any  
 connected component of $E_h$  to $M$ has the form  $T^n \times [s_1, s_2]$, where $0 
 \le s_2-s_1 <1$, for any $k\ne 0$ the mapping $\phi^k =\underbrace{\phi\circ ... \circ \phi}_{k \textrm{ times}}$ sends any   connected component of $E_h$ to another connected  component of $E_h$ so that 
  the procedure of taking quotient identifies  different connected  components of $E_h$ but 
   does not affect the dynamic on the components or  the topology of the components.  
  \end{remark} 
  \begin{remark} 
  The corrected version of \cite[Theorem 1.1(c)]{liu} could be: for every $h$ such that $V_{min} \le h \le V_{max}$, the topological entropy of the restriction of the Hamiltonian system to any connected component of  $E_h\subset \widetilde M$  is zero. For every $h>V_{max}$ the restriction of the Hamiltonian system to any connected component of  $E_h$  is positive\footnote{As we already mentioned above,  additional assumptions  on $A$ are assumed
  in \cite{liu}; these assumptions imply that the entropy is indeed positive for $h>V_{max}$ (which is not the case for example if $A= {\bf 0}$)}.  
  For $V_{min} \le h  <  V_{max}$ this statement follows from Theorem 1. For $h>V_{max}$, the statement was (correctly) proved in \cite{liu}; in this case, the topology of $E_h$ (the assumptions on $A$  assumed in \cite{liu}  imply  the exponential growth  of the fundamental group of $\widetilde M$)  forbids zero entropy. The remaining case $h=V_{max}$ can be obtained similar to the proof of Theorem 1, we explain  it in Remark 3.  \end{remark}

{\bf Proof of Theorem 1.}  The proof is similar to  (actually, is easier than)  the proof of  \cite[Theorem 1]{kruglikov}.  In order to prove that the topological entropy $h_{top}$ vanishes, we use the variational principle (see, for example, Theorem 4.5.3 of \cite{KH}):
 
 $$
h_{top}=\sup_{\mu\in \mathfrak{B}}h_\mu.
 $$
Here $\mathfrak{B}$ is the set of all invariant ergodic
probability measures on $E_h$  and $h_\mu$ is
the entropy of an invariant measure $\mu$. Recall that a measure
is called {\it ergodic}, if $\mu(B)(1-\mu(B))=0$ for all
$\mu$-measurable invariant Borel sets $B$. 

Therefore, in order to prove Theorem~1, it is sufficient 
 to prove that $h_\mu=0$ for all $\mu\in \mathfrak{B}$.
Fix one such measure and let $\Supp(\mu)$ be its support (the set
of $x\in M^{n+1}$ such that every  neighborhood
$U_\epsilon(x)$ has  positive measure).

We will use  that our  Hamiltonian  system is Liouville-integrable: the n+1 functions $p_1,...,p_n, H$ are integrals in the involution;  as it will be clear from the formulas \eqref{below} below, their  differentials  are linearly independent at almost every point of $T^*M$. 

Since the measure is ergodic, its support lies on a level surface
of every invariant continuous function. Then,
$\op{Supp}(\mu)$ is included into a Liouville leaf $\Upsilon$
(Recall that {\it a Liouville leaf}  is  a connected component of
the set $\{(x,p)\mid p_1=c_1,\dots,p_n=c_n, H= h\}$,  where $c_1,...,c_n$
are  constants.)

We say that a point $\xi \in T^*M$ is singular, if the differential of the integrals $p_1,...,p_n, H$ are linearly independent at this point.

Suppose at a  point $\xi\in \op{Supp}(\mu)$   is not singular. 
  Then, a small
neighborhood $U(\xi)$ of $\xi$ in $\op{Supp}(\mu)$
\begin{itemize}
\item has positive measure in $\mu$,
\item contains only nonsingular points.
\end{itemize}
We will show that these two conditions  imply that the entropy of
$\mu$ is zero.

In order to  do this, we consider the Poisson action of the the group $(\R^{n+1},+)$ on
$TM^{n+1}$: an  element $(a_1,...,a_{n+1})\in \R^{n+1} $ acts by time-one  shift
along the Hamiltonian vector field of the function 
$a_1p_1+...+a_np_n+ a_{n+1} H$.  Since the functions are  commuting integrals, the action is well-defined, smooth, symplectic,
 preserves the Hamiltonian of the geodesic flow, see   \cite[\S 49]{A} for details.

  By implicit function Theorem,  $\Upsilon$ is $n+1$-dimensional near
$\xi$. Denote  by $O(\xi)$ the orbit of the Poisson action of $(\R^{n+1},+)$
containing $\xi$. Since it is also $n+1$-dimensional, in a small
neighborhood of $\xi$ it coincides with $\Upsilon$. Thus,
$U(\xi)\subset O(\xi)$.

The orbits of the Poisson action and the dynamic on them are
well-studied  (see, for example,   \cite[\S 49]{A}).  There exists
a diffeomorphism to $$ T^k\times
\mathbb{R}^{n-k+1}=\underbrace{S^1\times ...\times S^1}_{k}\times
\underbrace{ \mathbb{R}\times...\times \mathbb{R}}_{n-k+1}$$ with
the standard coordinates $\phi_1,...,\phi_k\in (\mathbb{R}\ \
\textrm{mod} \ \ 2\pi)$, $t_{k+1},...,t_{n+1}\in \mathbb{R} $ such
that in these coordinates (the push-forward of) every trajectory
of the geodesic flow is given by the formula  \begin{equation} \label{below} \begin{array}{ll}
(\phi_1(\tau),...,\phi_k(\tau),&t_{k+1}(\tau),...,t_{n+1}(\tau))
                                            \\&=(\phi_1(0)+\omega_1\tau,\ ...\ ,\phi_k(0)+\omega_k\tau,
t_{k+1}(0)+\omega_{k+1}\tau,...,t_{n+1}(0)+\omega_{n+1}\tau),\end{array}\end{equation} where
 $\omega_1,...,\omega_{n+1}$ are constant on  $T^k\times
\mathbb{R}^{n-k}$.

We see that if at least one of the constants
$\omega_{k+1},...,\omega_{n+1}$ is not zero, every point  of $U(\xi)$  
 is wandering  in  $\Supp(\mu)$ 
 (see   \cite[\S3 in  Chapter 3]{KH} for definition), 
which contradicts the invariance of the measure.  
Then, the entropy of $\mu$ is zero.

If all constants $\omega_{k+1},...,\omega_{n+1}$ are zero, the coordinates $t_{k+1},...,t_{n+1}$ are constant on the trajectories of the geodesic flow. Since $\mu$ is ergodic, they are constant on the points of $\Supp(\mu)$. Then, $\Supp(\mu)$ 
is (diffeomorphic to) the torus $T^{\bar k}$ of dimension $\bar k\le k$, and 
the dynamics on  $\Supp(\mu)$  is (conjugate to)  the linear  flow  on  $T^{\bar k}$.  Then, the
entropy of $\mu$ is zero, see for example   \cite[Proposition~3.2.1]{KH}.

Suppose now that the point $\xi$ is singular, that is, the differentials of $p_1, ..., p_n, H$ are linearly dependent at $\xi$. In the coordinates $(p_1,...,p_{n+1}, x_1,...,x_{n+1})$,  the differentials of the   integrals are given by 

$$
\begin{array}{ccl}  dp_1& =&  (1,0,...,0, \underbrace{0,...,0}_{n+1})\\
                        &\vdots&    \\
                    dp_n&=&   (0,...,1,0, \underbrace{0,...,0}_{n+1})\\
                    dH &=&     (p_1,...,p_{n+1},\underbrace{0,...,0}_n, \tfrac{dV(x_{n+1})}{dx_{n+1}}).\end{array}
$$  
We see that the differentials are linearly dependent at $\xi=(p,x)$, if and  only if $p_{n+1}=0$ and $\tfrac{dV(x_{n+1})}{dx_{n+1}}=0$.  Take a trajectory $(p(t), x(t))$ passing through the point $\xi$. At the point $\xi$, in view of     $p_{n+1}=0$, we have $\tfrac{dx_{n+1}(t)}{dt}=0$. Since every point of the trajectory is singular,  at every point of the trajectory we have $p_{n+1}=0$. Then, $\tfrac{dx_{n+1}(t)}{dt}=0$ at every $t$  implying 
 $x_{n+1}=\const$ along the trajectory. Then,  the support of the measure $\mu$ lies on $E_h \cap \{(p,x)\in T^*M\mid x_{n+1}=\const\}$. 

   Now, since $\tfrac{dV(x_{n+1})}{dx_{n+1}}=0$ for $x_{n+1}=\const$, 
    the restriction of the Hamiltonian system to $E_h \cap \{(p,x)\in T^*M\mid x_{n+1}=\const\}$
     is actually (a subsystem of)  the geodesic flow of the flat metric 
   $$
\sum_{i,j=1}^nQ_{ij} dx_idx_j
$$ 
on the torus $T^n = \{x\in M \mid x_{n+1}=\const\}$.  
    Then, the entropy of $\mu$ is zero.  Theorem 1 is proved.  
\begin{remark} 
Let us now discuss   the case $h=V_{max}$ in the context of the corrected version (see Remark 2) 
of \cite[Theorem 1.1(c)]{liu}. At $(p,x)\in E_{V_{max}}$ 
 such that $V(x)= V_{max}$ we have  $p=0$.  
 Then, the  set   $E_{V_{max}} \cap  \{(p,x) \in T^*\widetilde M \mid  V(x)= V_{max}\}$ is  an invariant subset containing stable points only; in particular 
every invariant ergodic probability measure $\mu$ on this subset has zero entropy. 
Now, consider the set $\{ (p,x)\in E_{V_{max}} \mid V(x)<V_{max}\}.  $ The set is invariant as  the compliment to an invariant set. As we have shown in the proof of Theorem 1, 
every invariant ergodic probability measure $\mu$ on this subset has zero entropy. Finally, 
 the topological  entropy  of the whole system on $E_{V_{max}}\subset T^*\widetilde M$  is zero as we claimed.

\end{remark}

{\footnotesize \hspace{-10pt}\textsc{ {\sc Vladimir S. Matveev:} 
  Institute of Mathematics, \\
   07737 Jena Germany.
   {\tt vladimir.matveev@uni-jena.de}}}\\

\end{document}